\input amstex
\input Amstex-document.sty

\def\ZZ{{\Bbb Z}}

\def\CA{{\Cal{A}}}
\def\CG{{\Cal{G}}}

\def\ra{{\rightarrow}}

\def\CO{\Cal{O}}

\def\QQ{{\Bbb Q}}

\def\Qp{{\Bbb Q}_p}

\def\Qbar{{\overline{\Bbb Q}}}
\def\Qlb{{\overline{\Bbb Q}_{\ell}}}
\def\Fbar{{\overline{\Bbb F}_p}}

\def\ad{{\bold A}}

\def\CC{{\Bbb C}}
\def\PP{{\Bbb P}}

\define\isoarrow{{~\overset\sim\to\longrightarrow~}}
\def\-{{-1}}

\pageno 583

\def\shorttitle{On the Local Langlands Correspondence}
\def\shortauthor{Michael Harris}

\topmatter
\title\nofrills{\boldHuge On the Local Langlands Correspondence}
\endtitle

\author \Large Michael Harris* \endauthor

\thanks *Institut de Math\'ematiques de Jussieu-UMR CNRS 7586, Universit\'e Paris
7.  Membre, Institut Universitaire de France, France. E-mail: harris\@math.jussieu.fr
\endthanks

\abstract\nofrills \centerline{\boldnormal Abstract}

\vskip 4.5mm

{\ninepoint The local Langlands correspondence for $GL(n)$ of a non-Archimedean local field $F$ parametrizes
irreducible admissible representations of $GL(n,F)$ in terms of representations of the Weil-Deligne group $WD_F$
of $F$.  The correspondence, whose existence for $p$-adic fields was proved in joint work of the author with R.
Taylor, and then more simply by G. Henniart, is characterized by its preservation of salient properties of the two
classes of representations.

The article reviews the strategies of the two proofs. Both the author's proof with Taylor and Henniart's proof are
global and rely ultimately on an understanding of the $\ell$-adic cohomology of a family of Shimura varieties
closely related to $GL(n)$.  The author's proof with Taylor provides models of the correspondence in the
cohomology of deformation spaces, introduced by Drinfeld, of certain $p$-divisible groups with level structure.

The general local Langlands correspondence replaces $GL(n,F)$ by an arbitrary reductive group $G$ over $F$, whose
representations are conjecturally grouped in packets parametrized by homomorphisms from $WD_F$ to the Langlands
dual group ${}^LG$. The article describes partial results in this direction for certain classical groups $G$, due
to Jiang-Soudry and Fargues.

The bulk of the article is devoted to motivating problems that remain open even for $GL(n)$.  Foremost among them
is the search for a purely local proof of the correspondence, especially the relation between the Galois-theoretic
parametrization of representations of $GL(n,F)$ and the group-theoretic parametrization in terms of
Bushnell-Kutzko types. Other open questions include the fine structure of the cohomological realization of the
local Langlands correspondence:  does the modular local Langlands correspondence of Vigneras admit a cohomological
realization?

\vskip 4.5mm

}

\endabstract
\endtopmatter

\document

\baselineskip 4.5mm \parindent 8mm

\specialhead \noindent \boldLARGE Introduction \endspecialhead

Compared to the absolute Galois group of a number field, e.g.
$Gal(\Qbar/\QQ)$, the Galois group $\Gamma_F$ of a non-archimedean
local field $F$ has a ridiculously simple structure.  Modulo  the
inertia group $I_F$, there is a natural isomorphism
$$\Gamma_F/I_F \isoarrow Gal(\bar{k}_F/k_F),$$
where $k_F$ is the residue field of $F$.
Then $Gal(\bar{k}_F/k_F)$
is topologically generated by the geometric Frobenius
$Frob(x) = x^{\frac{1}{q}}$, where $q = |k_F| = p^f$ for $p$ prime.
The inertia group has a two step filtration,
$$1 \rightarrow P_F \rightarrow I_F \rightarrow \prod_{\ell \neq p}
\ZZ_{\ell} \ra 1,$$
where the wild ramification group $P_F$ is a pro-$p$ group.

Thus if $\sigma: \Gamma_F \ra GL(n,\CC)$ is a continuous
homomorphism, $n \geq 1$, then the image of $\sigma$ is solvable,
and $\sigma(P_F)$ is nilpotent.  This is still true when $\sigma$ is
a finite-dimensional complex representation of the {\it Weil group},
the subgroup
$W_F \subset \Gamma_F$ of elements whose image in $Gal(\bar{k}_F/k_F)$ is
an integral power of $Frob$.  Despite this simplicity, our understanding
of the set of equivalence classes of $n$-dimensional
representations of $W_F$ is far from complete, at least when $p$ divides $n$.

The reciprocity map of local class field theory:
$$F^{\times} \isoarrow W_F^{ab},$$
identifies the set $\CG(1,F)$ of one-dimensional representations
of $W_F$ with the set $\CA(1,F)$ of irreducible representations of
$F^{\times} = GL(1,F)$.  More than a simple bijection, this identification
respects a number of salient structures, and its behavior with respect
to field extensions $F'/F$ is well understood.  Moreover, it is compatible,
in a straightforward way, with global class field theory, and
was historically first derived as a consequence of the latter.

A simple special case of Langlands' functoriality principle is the
so-called strong Artin conjecture, which identifies the Artin L-function
attached to an irreducible $n$-dimensional representation of $Gal(\Qbar/K)$,
for a number field $K$, as the $L$-function of a cuspidal automorphic
representation of $GL(n)_K$.  As a local counterpart, Langlands
proposed a parametrization of irreducible admissible
representations of reductive groups over the local field $F$ in terms
of representations of
$W_F$.  The prototypical example is the local Langlands conjecture for $GL(n)$.
By analogy with the case $n = 1$, the set of equivalence classes of
irreducible admissible representations of $G_n = GL(n,F)$ is denoted
$\CA(n,F)$.  By
$\CG(n,F)$ we denote the set of equivalence classes of $n$-dimensional
representations, not of $W_F$ but rather of the Weil-Deligne group $WD_F$,
and only consider representations for which any lifting of $Frob$
acts semisimply.
Then the general local Langlands conjecture for $GL(n)$, in its crudest form,
asserts the existence of a family of bijections, as $F$ and $n$ vary:
$$\sigma = \sigma_{n,F}: \CA(n,F) \isoarrow \CG(n,F).  \tag 0.1$$
A normalization condition is that the central character $\xi_{\pi}$
of $\pi \in \CA(n,F)$
correspond to $\det \sigma(\pi)$ via local class field theory.

The first general result of this type was proved by Henniart [He1].
Early work of Bernstein and Zelevinsky reduced (0.1) to the existence
of bijections
$$\sigma = \sigma_{n,F}: \CA_0(n,F) \isoarrow \CG_0(n,F), \tag 0.2$$
where $\CG_0(n,F)$ are the irreducible representations of $W_F$ and
$\CA_0(n,F)$ is the supercuspidal subset of $\CA(n,F)$.  Both sides
of (0.2) are homogeneous spaces under $\CA(1,F)$, and thus
under its subset $\CA^{unr}(1,F)$ of unramified characters $\chi$ of
$F^{\times}$: if $\pi \in \CA_0(n,F)$ (resp. $\sigma \in
\CG_0(n,F)$), we denote by
$\pi\otimes \chi$ (resp. $\sigma\otimes \chi$)
the tensor product of $\pi$ (resp. $\sigma$) with the one-dimensional
representation
$\chi\circ \det$ of $G_n$ (resp. with the character $\sigma_{1,F}(\chi)$
of $W_F$).  Each $\CA^{unr}(1,F)$-orbit
on either side of (0.2) has a discrete invariant, the Artin conductor
$a(\pi)$, resp. $a(\sigma)$,
and the sets $\CA_0(n,F)[a]$, resp. $\CG_0(n,F)[a]$ of orbits with given
Artin conductor $a$ are known to be finite.  The main theorem of
[He1] is the  {\it numerical local Langlands correspondence}
$$|\CA_0(n,F)[a]| = |\CG_0(n,F)[a]|, \tag 0.3$$
established by painstakingly
counting both sides.

It has been known for some time that a family of bijections (0.2),
compatible with
Artin conductors and twists by $\CA(1,F)$, is not unique.  Henniart
showed (the Uniqueness Theorem, [He2]) that at most one normalized bijection
is compatible with contragredients and twists and satisfies the condition:
$$L(s,\pi\otimes\pi') = L(s,\sigma(\pi)\otimes\sigma(\pi')); \quad\quad
\varepsilon(s,\pi\otimes\pi',\psi) =
\varepsilon(s,\sigma(\pi)\otimes\sigma(\pi'),\psi) \tag 0.4$$
for $\pi \in \CA_0(n,F)$, $\pi' \in \CA_0(n',F)$, $n' < n$.  Here
$\psi: F \ra \CC^{\times}$ is
a non-trivial character.  The $L$- and $\varepsilon$-factors are
defined on the automorphic
side in [JPS, Sh]; on the Galois side by Langlands and Deligne (cf.
[D]).  It is in this version
that the local Langlands conjecture for $GL(n)$ has finally been established:
for fields of positive characteristic in [LRS], and for $p$-adic fields in
[HT], followed shortly thereafter by [He3] (see also [He5]).

\specialhead \noindent \boldLARGE 1.  Compatibility with global correspondences \endspecialhead

As in the first proofs of local class field theory, the bijections (0.2) are
constructed in [LRS,HT,He3] by local specialization of maps for certain
global fields $E$:
$$\sigma = \sigma_{n,E}: \CA^{good}(n,E) \hookrightarrow \CG(n,E). \tag 1.1$$
Here $E$ is supposed to have a place $w$ such that $E_w \isoarrow F$,
$\CA^{good}(n,E)$ is a class of cuspidal automorphic representations
of $GL(n)_E$ chosen to fit the circumstances, and $\CG(n,E)$ can
be taken to be the set of equivalence classes of compatible families
of $n$-dimensional semi-simple $\lambda$-adic representations of
$Gal(\bar{E}/E)$.
In particular, both sides of (1.1) as well as (0.2) are taken with
$\ell$-adic, rather than complex, coefficients; this does not change the
problem in an essential way.

The map $\sigma$ of (0.1) is particularly simple for unramified
representations.  An unramified $\tau \in \CG(n,F)$ is given by an
unordered $n$-tuple $(\chi_1,\dots,\chi_n)$ of unramified characters of
$W_F^{ab} \isoarrow F^{\times}$.  Ordering the $\chi_i$ arbitrarily, we obtain
a character $\chi$ of the Levi subgroup $G_1^n$ of a Borel subgroup
$B \subset G_n$.
The element of $\CA(n,F)$ corresponding to $\tau$ is then
the unique subquotient $\pi(\tau) = \sigma^{-1}(\tau)$
of the normalized induced representation $Ind_B^{GL(n,F)} \chi$
containing a vector fixed by $GL(n,\Cal{O}_F)$, $\Cal{O}_F$ the
integer ring of $F$.
This defines a bijection, a special case of
the {\it Satake parametrization}, between the unramified subset
$\CG^{unr}(n,F)$
and the unramified (spherical) representations $\CA^{unr}(n,F)$ of $G_n$.

Fix an automorphic representation $\Pi = \otimes_v \Pi_v$ of $GL(n)_E$.
The representation $\sigma_{n,E}(\Pi)$, when it exists, should
have the property that
$$\sigma_{n,E}(\Pi)~|~_{W_{E_v}} = \sigma_{n,E_v}(\Pi_v)  \tag 1.2$$
for almost all $v$ such that $\Pi_v \in \CA^{unr}(n,E_v)$; i.e.,
all but finitely many $v$.  By Chebotarev density, this
determines $\sigma_{n,E}(\Pi)$ uniquely.  One can then hope that
\proclaim{Hope 1.3} $\sigma_{n,E}(\Pi)_{W_{E_v}}$ depends only on
$F$ and $\Pi_v$
  for {\rom all} $v$,
\endproclaim
\noindent including $v = w$, the place of interest.  Setting
$\sigma_{n,F}(\Pi_v) = \sigma_{n,E}(\Pi)_{W_{E_v}}$, one then
needs to show that \proclaim{1.4}  For any $\pi \in \CA_0(n,F)$
there exists $\Pi \in \CA^{good}(n,E)$, for some $E$, with $\Pi_w
\simeq \pi$;
\endproclaim
\proclaim{1.5}  For $\Pi \in \CA^{good}(n,E)$, $\Pi' \in
\CA^{good}(n',E)$, the completed $L$-function \linebreak
$\Lambda(s,\sigma_{n,E}(\Pi)\otimes \sigma_{n',E}(\Pi'))$
satisfies the functional equation
$$\Lambda(s,\sigma_{n,E}(\Pi)\otimes \sigma_{n',E}(\Pi'))
= \varepsilon(s,\sigma_{n,E}(\Pi)\otimes
\sigma_{n',E}(\Pi'))\Lambda(1-s,\check{\sigma}_{n,E}(\Pi)\otimes
\check{\sigma}_{n',E}(\Pi'));$$
$$\varepsilon(s,\sigma_{n,E}(\Pi)\otimes \sigma_{n',E}(\Pi')) =
\prod_v \varepsilon_v(s,\sigma_{n,E}(\Pi)\otimes
\sigma_{n',E}(\Pi'),\psi_v)$$
is the product of local Deligne-Langlands $\varepsilon$ factors.
\endproclaim
Here $\check{}$ denotes contragredient.  The local additive
characters $\psi_v$ are assumed to be the local components of a
continuous character of $\ad_E/E$. \proclaim{1.6}  The map
$\sigma = \sigma_{n,F}: \CA_0(n,F) \ra \CG(n,F)$ \roster
\item"(i)" takes values in $\CG_0(n,F)$;
\item"(ii)" defines a bijection $\CA_0(n,F) \leftrightarrow \CG_0(n,F)$;
\item"(iii)" satisfies the remaining requirements of a local
Langlands correspondence, especially (0.4).
\endroster
\endproclaim

The main burden of [LRS] is the construction of a class
$\CA^{good}(n,E)$ large enough to
satisfy (1.4): now a moot point, since Lafforgue has proved that all cuspidal
automorphic representations of $GL(n)$ of a function field are
``good" in this sense.
The  $\CA^{good}(n,E)$ in [LRS] are the automorphic representations that
contribute to the cohomology of an appropriate Drinfeld modular
variety, constructed
from scratch for the occasion, attached to the multiplicative group
of a division
algebra of dimension $n^2$ over $E$, unramified at the chosen $w$.
Property (1.5)
in this case follows from general results of Deligne in [D], valid
only in equal characteristic.
Now by (1.2), for a sufficiently large set $S$ of places of $E$ we have
$$\prod_{v \notin S} L(s,\Pi_v \times \Pi'_v)  = \prod_{v \notin S}
L(s,\sigma_{n,E_v}(\Pi_v)\otimes \sigma_{n',E_v}(\Pi'_v)), \tag 1.7$$
where the left-hand side is the Rankin-Selberg L-function.
Completing the latter to $\Lambda(s,\Pi\otimes \Pi')$ and applying
[JPS] or [Sh],
we find the functional equation
$$\Lambda(s,\Pi\otimes \Pi') = \prod_v \varepsilon_v(s,\Pi\otimes
\Pi',\psi_v)\Lambda(1-s,\check{\Pi}\otimes \check{\Pi}'). \tag 1.8$$
In other words, the partial $L$-functions, identified via (1.8),
satisfy {\it two} functional equations (1.5) and (1.8).
An argument first used by Henniart then yields (0.4), and then (1.3)
and the full local Langlands conjecture
follow from the Uniqueness Theorem of [He2].

When $F$ is $p$-adic a class $\CA^{CK}(n,E)$ satisfying (1.2) is
implicitly defined by work of Clozel and Kottwitz [K,Cl1],
provided $E$ is a CM field.  For $\CA^{CK}(n,E)$ one can take
cuspidal automorphic representations $\Pi$, cohomological at all
archimedean primes, square integrable at several finite primes
other than $w$, and such that $\check{\Pi} \simeq \Pi^c$, where
$c$ denotes conjugation of $E$ over its maximal totally real
subfield.  However, the Galois-theoretic functional equation
(1.5) is only available a priori when $\sigma_{n,E}(\Pi)$ is
associated to a global complex representation of the Weil group
of $E$; i.e. when $\sigma_{n,E}(\Pi)$ becomes abelian over a
finite extension of $E$.  The article [H2] showed that there were
enough $\Pi$ of this type in $\CA^{CK}(n,E)$. Denoting by
$\CA^{good}(n,E)$ the set of such $\Pi$, we find that (1.4) is
impossible as soon as $p$ divides $n$; however, an argument in
[H2], based on Brauer's theorem on induced characters and (0.3),
shows that (1.4) is true ``virtually," in the set of formal sums
with integral coefficients of elements of $\CA^{good}(n,E)$ for
varying $n$.  It then suffices to prove the following weak form
of (1.3), which occupies the bulk of [HT]:
\proclaim{Theorem 1.9
[HT]} For all $\Pi \in \CA^{CK}(n,E)$, the semisimplification
\linebreak $\sigma_{n,E}(\Pi)_{W_{E_v},ss}$ of
$\sigma_{n,E}(\Pi)_{W_{E_v}}$ depends only on $F$ and $\Pi_v$ for
{\rom all} $v$.
\endproclaim

More precisely, [HT] proves that
$\sigma_{n,E}(\Pi)_{W_{F},ss}$ can be calculated explicitly in the
vanishing cycles of
certain formal deformation spaces $M^{h}_{LT,F}$ defined by Drinfeld
(see \S 2).
Following [K,Cl1], the representations $\sigma_{n,E}(\Pi)$ are
initially realized in the cohomology of certain
Shimura varieties with canonical models over $E$, and (1.9) is proved
by a study of their bad reduction at $w$.
Henniart soon realized that, for $\Pi \in \CA^{good}(n,E)$, the
purely local nature of $\sigma_{n,E}(\Pi)_{W_{E_v},ss}$,
and hence the definition of a map $\sigma_{n,F}$, could be derived
directly from (1.5) and
from the results of [He1,He2].  Though [He3] dispenses with the
geometry, it is still a global argument inasmuch
as it relies on [H2], which in turn depends on [K,Cl1] and [CL].
the conditional base change results of [CL].

A global consequence of Theorem 1.9 is the {\it Generalized Ramanujan
Conjecture} for the
automorphic representations in $\CA^{CK}(n,E)$:  if $\Pi \in
\CA^{CK}(n,E)$ and is unitary then
its local component $\Pi_v$ is tempered at every finite prime $v$.
Clozel in [Cl1] already showed
this to be true for almost all unramified $v$.
Generalizing a method developed
by Lubotzky, Phillips, and Sarnak for the $2$-sphere, Clozel [Cl2]
uses the version of the Generalized Ramanujan Conjecture
proved in [HT] to obtain effective constructions of families of
equidistributed points on odd-dimensional spheres.

With (0.1) out of the way, we can propose the following improvement of (1.3):

\proclaim{Problem 1}  {\rm Show that}
$$\sigma_{n,E}(\Pi)_{W_{E_v}} \isoarrow \sigma_{n,E_v}(\Pi_v). \tag 1.10$$
\endproclaim

For $n = 2$ this was established by Carayol
assuming standard conjectures
on the semisimplicity of Frobenius.  Theorem 1.9 shows that it holds
up to semisimplification.
The techniques of [HT], like the earlier work of Kottwitz treating
unramified places, is based on
a comparison of trace formulas, and cannot detect the difference
between two representations with
the same semisimplification.  Assuming semisimplicity of Frobenius,
the equality (1.10) follows easily from Theorem 1.9
and Deligne's conjecture, apparently inaccessible, on the
purity of the monodromy weight filtration.

\smallskip

\subhead{Compatibility with functoriality}
\endsubhead  Given cuspidal automorphic representations $\Pi_i$ of
$GL(n_i)_E$, for $i = 1, 2, \dots, r$, and a
homomorphism $\rho: GL(n_1)\times\dots \times GL(n_r) \ra GL(N)$ of
algebraic groups, Langlands functoriality predicts
the existence of an automorphic representation
$\rho_*(\Pi_1\otimes\dots \Pi_r)$, not necessarily
cuspidal, of $GL(N)_E$, such that, for almost all unramified places $v$ of $E$,
$$\sigma_{N,E_v}(\rho_*(\Pi_1\otimes\dots \Pi_r)_v) =
\rho\circ(\otimes_{i = 1}^r \sigma_{n_i,E_v}(\Pi_{i,v})). \tag 1.11$$
In recent years this has been proved
for general number fields $E$ in several important special cases:
the tensor products
$GL(2)\times GL(2) \ra GL(4)$ (Ramakrishnan) and $GL(2)\times GL(3)
\ra GL(6)$ (Kim-Shahidi), and the
symmetric powers $Sym^3: GL(2) \ra GL(4)$ (Kim-Shahidi) and $Sym^4:
GL(2) \ra GL(5)$  (Kim).
It has been verified in all four cases that (1.11) holds for all $v$.

\subhead{Construction of supercuspidal representations by ``backwards lifting"}
\endsubhead

The unitary representation $\pi \in \CA_0(n,F)$ is isomorphic to its
contragredient
if and only if the local factor $L(s,\pi\times\pi)$ has a pole at $s
= 0$, which
is necessarily simple.  The local
factor can be decomposed as a product:
$$L(s,\pi\times \pi) = L(s,\pi,Sym^2)L(s,\pi,\wedge^2), \tag 1.12$$
where the two terms on the right are defined for unramified $\pi$ by
Langlands and in general by Shahidi.  Only one of the factors on the right
has a pole. Using the class $\CA^{good}(n,E)$ of automorphic
representations, Henniart has shown that it is the first factor
(resp. the second
factor) if and only if $\sigma(\pi)$ is orthogonal (resp.
symplectic); the symplectic
case only arises for $n$ even.
One thus expects that $\pi$ is obtained by functorial transfer from
an $L$-packet
of a classical group $G$ over $F$, via the map
of $L$-groups ${}^LG \ra GL(n,\CC)$, where ${}^LG = SO(n,\CC)$, resp.
$Sp(n,\CC)$,
if the first, resp. the second factor in (1.12) has a pole at $s = 0$.

In particular, when $n = 2m$ and $L(s,\pi,\wedge^2)$ has a pole at $s
= 0$, $\pi$
should come from an $L$-packet on the split group $SO(2m+1,F)$.
Using a local analogue of the
method of ``backwards lifting," or automorphic descent, due to
Ginzburg, Rallis, and Soudry, Jiang and Soudry have constructed a generic
supercuspidal representation $\pi'$ of $SO(2m+1,F)$ for every
$\pi \in \CA_0(n,F)$ with the indicated pole.  More generally, they have
obtained a complete parametrization of generic representations of
split $G = SO(2m+1,F)$
in terms of Langlands parameters $WD_F \ra {}^LG$ [JS].  These results should
certainly extend to other classical groups.

\specialhead \noindent \boldLARGE 2.  Cohomological realizations of the local  correspondence \endspecialhead

The theory of the new vector implies easily
that any irreducible admissible representation $\pi \in \CA(n,F)$ has
a rational model
over
the field of definition of its isomorphism class:
the Brauer obstruction is trivial for $G_n$.
The analogous assertion fails for representations in $\CG(n,F)$.
Thus one cannot expect the existence
of a space $\Cal{M}$, with a natural action of $G_n\times W_F$, whose
cohomology of whatever sort realizes the local Langlands
correspondence, as an identity of virtual representations
$$\sigma_{n,F}(\pi) = \pm[Hom_{G_n}(H_c(\Cal{M}),\pi)] := \pm \sum_i
(-1)^i Hom_{G_n}(H^i_c(\Cal{M}),\pi). \tag 2.1$$

We add a third group to the picture by taking $J$ to be an inner form of $G_n$,
the multiplicative group of a central simple algebra $D$ over $F$ of
dimension $n^2$, with Hasse invariant $\frac{r_D}{n}$.  The set
$\CA(n,F)$ contains a subset $\CA_{(2)}(n,F)$
of discrete series representations, character twists
of those realized in the regular representation on $L_2(G_n)$ (modulo
center).  The set
$\CA(J)$ of equivalence classes of irreducible admissible representations
contains an analogous subset $\CA_{(2)}(J)$, equal to $\CA(J)$ if
$D$ is a division algebra.  The {\it Jacquet-Langlands correspondence }
[R,DKV] is a bijection
$JL: \CA_{(2)}(G_n) \isoarrow  \CA_{(2)}(J)$
determined by the identity of distribution characters
$$\chi_{\pi}(g) = \varepsilon(J) \chi_{JL(\pi)}(j), ~~\pi \in
\CA_{(2)}(G) \tag 2.2$$
if $\varepsilon(J) = \pm 1$ is the Kottwitz sign and $g$ and $j$ are elliptic
regular elements with the same eigenvalues.  When $r_D = 1$ there are
two spaces $\tilde{\Omega}^n_F$
and $\Cal{M}^n_{LT,F}$ with natural $G_n\times J$-actions.  The
former is a countable
union, indexed by $\ZZ$, of copies of the profinite \'etale
cover $\tilde{\Omega}^{n,0}_F$ of the rigid-analytic upper half space
$\Omega^n_F = \PP^{n-1}(\CC_p) - \PP^{n-1}(F)$, defined
by Drinfeld in [D2].  The latter is the rigid generic fiber of the
formal deformation space $M^n_{LT,F}$ of the one-dimensional height
$n$ formal
$\CO_F$-module with Drinfeld level structures of all degrees [D1].  A
relation analogous to (2.1)
was conjectured by Carayol in [C1], with $\pm = (-1)^{n-1}$:
\proclaim{Theorem 2.3}  For $\pi$ supercuspidal
$$\aligned &\sigma^{\#}(\pi)\otimes JL(\pi) =
\pm[Hom_{G_n}(H_c(\tilde{\Omega}^n_F),\pi)] \\
&\sigma^{\#}(\pi)\otimes \pi =
\pm[Hom_{J}(H_c(\Cal{M}^n_{LT,F}),JL(\pi))].\endaligned $$
\endproclaim

The notation $\sigma^{\#}(\pi)$ indicates that $\sigma(\pi)$ has been
twisted by an elementary factor.  We  use the rigid-analytic
\'etale cohomology introduced by Berkovich [B] with coefficients
in $\Qlb$, $\ell \neq p$.  For $\Cal{M}^n_{LT,F}$ this can be interpreted
as a space of vanishing cycles for the formal deformation
space, viewed as a formal scheme over $Spf(\CO_F)$.
The case of $\tilde{\Omega}^n_F$ was proved in [H1], using the
existence of Shimura varieties
admitting rigid-analytic uniformizations by $\tilde{\Omega}^n_F$.
This has recently been extended to $F$ of equal characteristic
by Hausberger [Hau].  The case of $\Cal{M}^n_{LT,F}$, again for $\pi$
supercuspidal, was initially treated by Boyer [Bo] in the
equal-characteristic case.  The analogous statement for $F$ $p$-adic,
and for any $\pi$, is the logical starting point
of the proof of Theorem 1.9 in [HT].

Theorem 2.3 is extended in [HT] to general $\pi \in \CA_{(2)}(G)$.
The explicit formula
for the alternating sum of the
$Hom_{J}(H^i_c(\Cal{M}^n_{LT,F}),JL(\pi))$ is awkward
but yields a simple expression for
$$\multline \sum_{i,j} (-1)^{i+j}
Ext^j_G(Hom_{J}(H^i_c(\Cal{M}^n_{LT,F}),JL(\pi)),\pi)
\\ =  \sum_{i,j,k} (-1)^{i+j+k}
Ext^j_G(Ext^k_{J}(H^i_c(\Cal{M}^n_{LT,F}),JL(\pi)),\pi)
\endmultline\tag 2.4$$
in terms of the semisimplification of $\sigma(\pi)$.
An analogous {\it conjectural} expression for {\it individual}
$H^i_c(\tilde{\Omega}_F)$
has been circulating for several
years and should appear in a forthcoming joint paper with Labesse.
Faltings has proved [F2]  that the spaces $\tilde{\Omega}_F$ and
$\Cal{M}^n_{LT,F}$ become isomorphic after $p$-adic
completion of the latter.  Thus the two questions in the following
problem reduce to a single question:

\proclaim{Problem 2}  {\rm Determine the individual
representations $H^i_c(\Cal{M})$, and the spaces\linebreak
$Ext^j_{G_n}(H^i_c(\tilde{\Omega}_F),\pi)$ and
$Hom_{J}(H^i_c(\Cal{M}^n_{LT,F}),JL(\pi))$ for all $i, j$,  all
$\pi \in \CA(n,F)$.  In particular, show that
$Ext^j_{G_n}(H^i_c(\tilde{\Omega}_F),\pi)$ vanishes unless there
exists $\pi' \in \CA_{(2)}(n,F)$ such that $\pi$ and $\pi'$
induce the same character of the Bernstein center.}
\endproclaim

The results of [HT] imply that, for any $\pi \in \CA_{(2)}(n,F)$,
with Bernstein
character $\beta_{\pi}$, the Bernstein
center acts on
$\sum_i (-1)^i Hom_{J}(H^i_c(\Cal{M}^n_{LT,F}),JL(\pi))$
via $\beta_{\pi}$.

For $\pi$ supercuspidal it is known in all cases that the spaces in
Problem 2 vanish
for $i \neq n-1$ (and for $j \neq 0$).  This vanishing property
should characterize supercuspidal $\pi$ among representations
in $\CA_{(2)}(n,F)$.  When $\pi$ is the Steinberg representation, the
$H^i_c(\Omega_F)$, as well
as the corresponding $Ext$ groups, are calculated explicitly in [SS].
The calculation in [SS] is purely local, whereas the vanishing
outside the middle
degree for $\pi \in \CA_0(n,F)$ is based on properties of automorphic forms.

\proclaim{Problem 3}  {\rm Find a purely local proof of the
vanishing property for $\pi \in \CA_0(n,F)$.}

\endproclaim

The covering group of $\tilde{\Omega}^{n,0}_F$ over $\Omega^n_F$ can
be identified with the
maximal compact subgroup $J^0 \subset J$.  Thus
$H^i_c(\tilde{\Omega}_F)$ can be written
as a sum $\oplus_{\tau} H^i_c(\tilde{\Omega}_F)[\tau]$ of its
$\tau$-isotypic components,
where $\tau$ runs over irreducible representations of $J^0$ or,
equivalently, over inertial
equivalence classes of representations of $J$.  Closely related to
Problem 3 is the following

\proclaim{Problem 4}  {\rm Characterize $\tau \in \CA(J)$ such
that $JL^{-1}(\tau) \in \CA_0(n,F)$.  Equivalently, calculate the
Jacquet functors of the $G_n$-spaces
$H^i_c(\tilde{\Omega}_F)[\tau]$ geometrically, in terms of
$\tau$.}
\endproclaim

When $n$ is prime $JL^{-1}(\tau) \in \CA_0(n,F)$ if and only if $\dim
\tau > 1$; when $\dim \tau = 1$
$JL^{-1}(\tau)$ is a twist of the Steinberg representation.  For general
$n$ practically nothing is known.

\subhead Results of L. Fargues [Fa]
\endsubhead  For certain classical $F$-groups $G$, Rapoport and Zink,
using the deformation
theory of $p$-divisible (Barsotti-Tate) groups, have defined pro-rigid analytic
spaces $\Cal{M}$ admitting continuous $G\times J \times W_E$ actions on their
$\ell$-adic cohomology, where $J$ is an inner form of $G$ and $E$,
the {\it reflex field}
of $\Cal{M}$, is a finite extension of $F$ [RZ].  In [R] Rapoport
proposes a conjectural formula, which
he attributes to Kottwitz, for the discrete series contribution to
the virtual $G\times J \times W_E$-module $[H(\Cal{M})] = \sum_{i}
(-1)^i H^i_c(\Cal{M},\Qlb)$.
The pairs $(G,J)$ considered in [RZ] include $(G_n,D^{\times})$ with
general Hasse
invariant $\frac{r_D}{n}$, $G = J = GU(n)$, the quasi-split unitary
similitude group
attached to the unramified quadratic extension of $F$, and the symplectic
similitude group $G = GSp(2n,F)$.

\proclaim{Theorem 2.5 (Fargues)}  Suppose $F/\Qp$ unramified, $(G,J)
= (G_n,D^{\times})$,
with $(r_D,n) = 1$.  For any $\pi \in \CA_0(n,F)$ we have
$$\sum_i (-1)^i Hom_{J}(H^i_c(\Cal{M},\Qlb),JL(\pi))_0 = \pm
\pi\otimes\wedge^{?}\sigma(\pi)$$
up to a simple twist.  Here the subscript ${}_0$ denotes the
$G$-supercuspidal part
and $\wedge^{?}$ is a certain tensor product of exterior powers of
$\sigma(\pi)$ with total weight $r_D$,
depending on auxiliary data defining $\Cal{M}$.
\endproclaim

This confirms the Kottwitz-Rapoport conjectures
in the case in question.  For $G = J = GU(3)$ Rogawski has defined a
local Langlands correspondence via base change
to $GL(3)$.  In that case the supercuspidal representations of $G$
are grouped into
$L$-packets.  Fargues's techniques apply to this case as well, and he
obtains a version of
the Kottwitz-Rapoport conjectures, more difficult to state than
Theorem 2.5 (higher $Ext$'s are
involved, and the formula is averaged over $L$-packets).\footnote{The
statement of
the general Kottwitz-Rapoport conjectures in [H3] for general discrete series
representations is based on a misreading
of Rapoport's  use of the term ``discrete $L$ parameter".  The
correct conjecture should
involve the analogue of the alternating sum on the right-hand side of
(2.4), with $JL(\pi)$
replaced by $\pi'$ in  the $L$-packet associated to $\pi$.}  More
generally, Fargues' methods apply
to classical groups attached to Shimura varieties, whenever the trace
formula is
known to be stable and functorial transfer from $G$ to $GL(n)$ has
been established.

In contrast to [HT], Fargues' methods are essentially rigid-analytic,
and make no use
of equivariant regular integral models of Shimura varieties in wildly
ramified level -- fortunately so,
since such models are not known to exist.  Heuristically, the
characters of the representations of
$G$ and $J$ on $[H(\Cal{M})]$ can be related by applying a Lefschetz
trace formula to $\ell$-adic cohomology
of the rigid space $\Cal{M}$.  This approach, which in principle provides
no information about the $W_F$ action, has been successfully applied
to $\tilde{\Omega}_F^n$
by Faltings in [F1], and to $M^n_{LT,F}$ by Strauch [S] when $n = 2$.
For higher $M^n_{LT,F}$,
and for the Rapoport-Zink spaces studied by Fargues, one does not yet
know how to deal with
wild boundary terms in Huber's Lefschetz formula [Hu] and its
higher-dimensional generalizations.

Using work of Oort and Zink on stratification of families of abelian
varieties and
the slope filtration for $p$-divisible groups, Mantovan [M] has
developed another
approach to the cohomology of Shimura varieties of PEL type.  Closer
in spirit to
[HT] than to [F], [M] obtains finer results on the geometry of the special
fiber and a description of the cohomology in ramified level similar to that
of [F].

\subhead{Cohomological realizations with torsion coefficients}
\endsubhead

It would be convenient if the following question
had an affirmative answer:

\proclaim{Question 5}  {\rm Is
$H^i_c(\tilde{\Omega}^n_F,\ZZ_{\ell})$ a torsion-free
$\ZZ_{\ell}$-module?}
\endproclaim

The global trace formula methods used in [H1] and [HT] to derive
Theorem 2.3 from an analysis of the cohomology of the ``simple"
Shimura varieties of the title of [K] are insensitive to torsion in
cohomology.  When $\ell > n$ it may be possible,
as in recent work of Mokrane and Tilouine,
to combine $\ell$-adic Hodge theory with the generalized
Eichler-Shimura congruence formula,
for the same ``simple" Shimura varieties, to answer Question 5.  For
$\ell \leq n$
completely new ideas are needed.

When $k$ is an algebraically closed field of characteristic $\ell
\neq p$, Vign\'eras
has defined a class of smooth supercuspidal representations
$\CA_{0,k}(n,F)$ of $G_n$
with coefficients in $k$, and has proved that they are in bijection with
the set $\CG_{0,k}(n,F)$ of irreducible $n$-dimensional
representations of $W_F$
over $k$ (see article in these Proceedings).  It is natural to expect that this
modular local Langlands correspondence is realized on the spaces
$H^{\bullet}_c(\Cal{M},k)$,
with $\Cal{M} = \tilde{\Omega}^n_F$ or $\Cal{M}^n_{LT,F}$.

\proclaim{Problem 6} {\rm Define a modular Jacquet-Langlands map
$\pi \mapsto JL(\pi)$ from $\CA_{0,k}(n,F)$ to
$k$-representations of $J$, and formulate the last sentence
precisely.  Does the virtual $W_F$-module
$$(-1)^{n-1}\sum_{i,j,k} (-1)^{i+j+k}
Ext^j_G(Ext^k_{J}(H^i_c(\Cal{M},k),JL(\pi)),\pi)$$ realize the
modular local Langlands correspondence?}
\endproclaim

Implicit in the second question is the assumption that the modular
Jacquet-Langlands
map can be extended to a wider class of $k$-representations of $G_n$, perhaps
including reduction $\pmod{\ell}$ of supercuspidal representations in
characteristic zero.
One can of course ask the same questions when $\ell = p$.  In this case we can
consider rigid (de Rham) cohomology, in the sense of Berthelot, as
well as $p$-adic \'etale
cohomology.  All three groups $G_n$, $J$, and $W_F$ have large
analytic families
of $p$-adic representations.  It is not at all clear whether the
$p$-adic cohomology of
$\tilde{\Omega}^n_F$ is sufficiently rich to account for all $p$-adic
deformations -- in
categories yet to be defined -- of a given representation occurring
in cohomology with
coefficients in $\Fbar$.
\bigskip

\specialhead \noindent \boldLARGE 3.  Explicit parametrization of supercuspidal representations \endspecialhead

\subhead{Distribution characters}
\endsubhead

The distribution character $\chi_{\pi}$, a locally integrable
function on the set of regular semisimple elements of $G_n =
GL(n,F)$, is the fundamental analytic invariant
of $\pi \in \CA(n,F)$.  For $\pi \in \CA_{(2)}(n,F)$, $\chi_{JL(\pi)}$, related
to $\chi_{\pi}$ by (2.2), extends continuously to an invariant
function on $J = D^{\times}$
provided $(r_D,n) = 1$, which we assume.  Under this hypothesis every
element of
$J$ is elliptic and every elliptic regular element $j$ is contained
in a unique maximal
torus $T(j)$, isomorphic to the multiplicative group of an extension
$K$ of $F$ of degree $n$.
Since $JL(\pi)$ is finite-dimensional, its restriction to $T(j)$
equals $\sum_{\xi} a_{\pi}(\xi)\xi$
where $\xi$ runs over characters of $K_j^{\times}$ and the
coefficients $a_{\pi}(\xi) = a_{\pi}(K,\xi)$ are
non-negative integers, almost all zero.  In this way $\pi \in
\CA_{(2)}(n,F)$ is determined by the integer-valued function
$a_{\pi}(K,\xi)$ where $K$ runs over degree $n$ extensions of $F$ and
$\xi$ over characters of $K^{\times}$.
Invariance entails the symmetry condition $a_{\pi}(K',^{\sigma}\xi) =
a_{\pi}(K,\xi)$
where $\sigma: K  \isoarrow K'$ is an isomorphism over $F$; in
particular, if $\sigma \in Aut_F(K)$.

\proclaim{Problem 7}  {\rm Express $a_{\pi}(K,\xi)$ in terms of
numerical invariants of $\sigma(\pi)$.}
\endproclaim

Of course $a_{\pi}(K,\xi) = 0$ unless $\xi|_{F^{\times}}$ coincides
with the central
character $\xi_{\pi}$ of $\pi$.  When $n = 2$ $a_{\pi}(K,\xi) \in
\{0,1\}$, and a theorem of Tunnell,
completed by H. Saito, relates the nonvanishing of $a_{\pi}(K,\xi)$
to the local constant
$\varepsilon(\frac{1}{2},\sigma(\pi)\otimes\xi^{-1},\psi)$.  For $n$
prime to $p$ a
conjecture of Reimann, following an earlier conjecture of
Moy, expresses $\chi_{\pi}$ in terms of $\sigma(\pi)$; work in
progress of Bushnell
and Henniart shows that this conjecture is almost right (probably up
to an unramified character
of degree at two).

\subhead{Parametrization via types}
\endsubhead

A fundamental theorem of Bushnell and Kutzko asserts
that every supercuspidal
$\pi$ can be obtained by compactly supported induction from a
finite-dimensional representation
$\tau$ of a subgroup $H \subset G_n$ which is compact modulo the
center $Z_n$ of $G_n$.
The pair $(H,\rho)$, called an {\it extended type}, is unique up to
conjugation by $G_n$.
The character $\chi_{\pi}$ can be obtained from $(H,\rho)$ by a
simple integral formula
[BH, (A.14)].

The outstanding open problem concerning the local Langlands
correspondence is undoubtedly

\proclaim{Problem 8} {\rm (a)  Define $\sigma(\pi)$ directly in
terms of $(H,\rho)$ (and vice versa).

(b)  Show directly that the definition of $\sigma$ in (a) has the
properties of a local Langlands correspondence.}
\endproclaim

Note that (b) presupposes a direct construction of the local
Galois constants.

Problem 8 formulates the hope, often expressed, for a purely local construction
of the local Langlands correspondence.  Bushnell, Henniart, and
Kutzko have made
considerable progress toward this goal.  Among other results, they
have obtained:

\roster
\item"$\bullet$"  A formula for the conductor $a(\pi\times\pi')$,
$\pi \in \CA_0(n,F)$, $\pi' \in \CA_0(n',F)$ [BHK];
\item"$\bullet$"  A purely local candidate for the base change map
$\CA(n,F) \ra \CA(n,K)$ when $K/F$ is a tame, not necessarily Galois
extension [BH, I], agreeing with Arthur-Clozel base change for $K/F$ cyclic;
\item "$\bullet$" A bijection between wildly ramified supercuspidal
representations of $G_{p^m}$ and wildly ramified\footnote{A wildly ramified
irreducible representation of $W_F$ is one that remains irreducible
upon restriction to $P_F$; a wildly ramified supercuspidal is one not
isomorphic to its twist by any non-trivial unramified character.}
representations
in $\CG_0(p^m,F)$, preserving local constants [BH, II].
\endroster
In each instance, the constructions and proofs are based primarily on
the theory of types.  A complete solution of Problem 8 remains elusive,
however, absent a better understanding of the local Galois
constants.

\proclaim{Question 9}  {\rm Can the types $(H,\rho)$ be realized
in the cohomology ($\ell$-adic or $p$-adic) of appropriate
analytic subspaces of $\tilde{\Omega}^n_F$ or $\Cal{M}^n_{LT,F}$?}
\endproclaim

Positive results for certain $(H,\rho)$ have been announced by
Genestier and Strauch, at least when $n = 2$.

\noindent {\bf Acknowledgments.} I thank R. Taylor, G. Henniart,
and L. Fargues for their comments on earlier versions of this
report.

\widestnumber \key {AAAA}

\specialhead \noindent \boldLARGE References \endspecialhead

More or less detailed accounts of the history of the local Langlands
conjecture,
and of its proofs, can already be found in the literature:  [Rd] and [Ku]
describe the problem and the work of Bernstein and Zelevinsky, while the
proofs are outlined in [C2,C3], [He4], [W], as well
as the introduction to [HT].

\ref \key AC \by  Arthur, J. and L. Clozel, \paper \rm {\it Simple algebras, base change, and the advanced theory
of the trace formula}, {\it Annals of Math. Studies}, {\bf 120}, Princeton:  Princeton University Press (1989)
\endref

\ref \key B \by  Berkovich, V.G. \paper \rm  \'Etale cohomology
for non-archimedean analytic spaces, {\it Publ. Math. I.H.E.S.},
{\bf 78}, 5--161 (1993) \endref

\ref \key Bo \by  Boyer, P. \paper \rm  Mauvaise r\'eduction de
vari\'et\'es de Drinfeld et correspondance de Langlands locale,
{\it Invent. Math.}, {\bf 138}, 573--629 (1999) \endref

\ref \key BHK \by  Bushnell, C., G. Henniart, and P. Kutzko
\paper \rm  Local Rankin-Selberg convolutions for $GL_n$:
explicit conductor formula, {\it J. Amer. Math. Soc.}, {\bf 11},
703--730 (1998) \endref

\ref \key BK \by  Bushnell, C. and P. Kutzko \paper \rm  The admissible dual of $GL(N)$ via compact open
subgroups, {\it Annals of Math. Studies}, {\bf 129} (1993) \endref

\ref \key BH \by  Bushnell, C. and G. Henniart \paper \rm  Local tame lifting for $GL(n)$, I.  {\it Publ. Math.
IHES}, {\bf } (1996); II:  wildly ramified supercuspidals, {\it Ast\'erisque}, {\bf 254} (1999) \endref

\ref \key C1 \by  Carayol, H. \paper \rm Non-abelian Lubin-Tate
theory,  in L. Clozel and J.S. Milne, eds., {\it Automorphic
Forms, Shimura varieties, and $L$-functions}, New York:  Academic
Press, vol II,  15--39 (1990) \endref

\ref \key C2 \by  Carayol, H. \paper \rm Vari\'et\'es de Drinfeld
compactes, d'apr\`es Laumon, Rapoport, et Stuhler, S\'eminaire
Bourbaki exp. 756 (1991-1992), {\it Ast\'erisque} {\bf 206}
(1992), 369--409 \endref

\ref \key C3 \by  Carayol, H. \paper \rm Preuve de la conjecture de Langlands locale pour $GL_n$:  Travaux de
Harris-Taylor et Henniart, S\'eminaire Bourbaki exp. 857. (1998-1999) \endref

\ref \key Cl1 \by  Clozel, L. \paper \rm Repr\'esentations
Galoisiennes associ\'ees aux repr\'esentations automorphes
autoduales de GL(n), {\it Publ. Math. I.H.E.S.}, {\bf 73},
97--145 (1991) \endref

\ref \key Cl2 \by Clozel, L. \paper \rm  Automorphic forms and the distribution of points on odd- \linebreak
dimensional spheres, (manuscript, 2001) \endref

\ref \key CL \by Clozel, L. and J.-P. Labesse \paper \rm  Changement de base pour les repr\'esentations
cohomologiques de certains groupes unitaires, appendix to J.-P. Labesse, {\it Cohomologie, stabilisation, et
changement de base}, {\it Ast\'erisque}, {\bf 257} (1999) \endref

\ref \key D \by  Deligne, P. \paper \rm  Les constantes des
\'equations fonctionnelles des fonctions $L$, {\it Modular
Functions of one variable II}, {\it Lect. Notes Math.}, {\bf
349}, 501--595 (1973) \endref

\ref \key DKV \by  Deligne, P., D.Kazhdan, and M.-F.Vigneras \paper \rm Repr\'{e}sentations des alg\`{e}bres
centrales simples $p$-adiques, in J.-N.Bernstein, P.Deligne, D.Kazhdan, M.-F.Vigneras, {\it Repr\'{e}sentations
des groupes r\'{e}ductifs sur un corps local}, Paris: Hermann (1984) \endref

\ref \key D1 \by  Drinfeld, V. \paper \rm  Elliptic modules, {\it
Math. USSR Sbornik}, {\bf 23}, 561--592 (1974)
\endref

\ref \key D2 \by  Drinfeld, V. \paper \rm  Coverings of $p$-adic
symmetric domains, {\it Fun. Anal. Appl.}, {\bf 10}, 107--115
(1976) \endref

\ref \key F1 \by Faltings, G. \paper \rm  The trace formula and
Drinfeld's upper halfplane, {\it Duke Math. J.}, {\bf 76},
467--481 (1994) \endref

\ref \key F2 \by Faltings, G. \paper \rm  A relation between two moduli spaces studied by V. G. Drinfeld,
preprint, 2001
\endref

\ref \key Fa \by Fargues, L. \paper \rm  Correspondances de Langlands locales dans la cohomologie des espaces de
Rapoport-Zink, th\`ese de doctorat, Universit\'e Paris 7 (2001) \endref

\ref \key Fu \by Fujiwara, K. \paper \rm Rigid geometry,
Lefschetz-Verdier trace formula and Deligne's conjecture, {\it
Invent. Math.}, {\bf 127}, 489--533 (1997) \endref

\ref \key H1 \by  Harris, M. \paper \rm Supercuspidal
representations in the cohomology of Drinfel'd upper half spaces;
elaboration of Carayol's program, {\it Invent. Math.}, {\bf 129},
75--119 (1997) \endref

\ref \key H2 \by  Harris, M. \paper \rm The local Langlands
conjecture for $GL(n)$ over a $p$-adic field, $n < p$, {\it
Invent. Math.}, {\bf 134}, 177--210 (1998) \endref

\ref \key H3 \by  Harris, M. \paper \rm Local Langlands correspondences and vanishing cycles on \linebreak Shimura
varieties, Proceedings of the European Congress of Mathematics, Barcelona, 2000. {\it Progress in Mathematics},
{\bf 201}, 407--427 (2001)
\endref

\ref \key HT \by  Harris, M. and R. Taylor \paper \rm  On the geometry and cohomology of some simple Shimura
varieties, {\it Annals of Math. Studies}, {\bf 161} (2002) \endref

\ref \key Hau \by  Hausberger, T.  Repr\'esentations cuspidales dans la cohomologie des \linebreak rev\^etements
de Drinfeld: preuve de la conjecture de Drinfeld-Carayol en \'egale caract\'eristique, preprint (2001) \endref

\ref \key He1 \by  Henniart, G. \paper \rm  La conjecture de
Langlands locale num\'erique pour $GL(n)$, {\it Ann. scient. Ec.
Norm. Sup}, {\bf 21}, 497--544 (1988) \endref

\ref \key He2 \by  Henniart, G. \paper \rm Caract\'erisation de
la correspondence de Langlands locale par les facteurs $\epsilon$
de paires, {\it Invent. Math.}, {\bf 113}, 339--350 (1993) \endref

\ref \key He3 \by  Henniart, G. \paper \rm Une preuve simple des conjectures de Langlands pour $GL(n)$ sur un
corps $p$-adique, {\it Invent. Math.}, (2000) \endref

\ref \key He4 \by Henniart, G. \paper \rm  A Report on the Proof
of the Langlands Conjectures for $GL(N)$ over $p$-adic Fields,
Current Developments in Mathematics 1999, International Press,
55--68 (1999) \endref

\ref \key He5 \by  Henniart, G. \paper \rm  Sur la conjecture de
Langlands locale pour $GL_n$, {\it J. Th\'eorie des Nombres de
Bordeaux}, {\it 13}, 167--187 (2001) \endref

\ref \key Hu \by  Huber, R. \paper \rm  Swan representations
associated with rigid analytic curves, {\it J. Reine Angew.
Math.}, {\bf 537}, 165--234 (2001) \endref

\ref \key JPS \by  Jacquet, H., I. I. Piatetski-Shapiro, and J.
Shalika \paper \rm  Rankin-Selberg convolutions, {\it Am. J.
Math.}, {\bf 105}, 367--483 (1983) \endref

\ref \key JS \by  Jiang, D. and D. Soudry \paper \rm  Generic representations and local Langlands reciprocity law
for $p$-adic $SO_{2n+1}$, preprint (2001) \endref

\ref \key K \by Kottwitz, R. \paper \rm On the $\lambda$-adic
representations associated to some simple Shimura varieties, {\it
Invent. Math.}, {\bf 108}, 653--665 (1992) \endref

\ref \key Ku \by Kudla, S. \paper \rm  The local Langlands
correspondence:  The non-archimedean case, {\it Proc. Symp. Pure
Math.}, {\bf 55}, part 2, 365--391 (1994) \endref

\ref \key LRS \by  Laumon, G., M. Rapoport, and U. Stuhler \paper
\rm $\Cal{D}$-elliptic sheaves and the Langlands correspondence,
{\it Invent. Math.}, {\bf 113}, 217--338 (1993) \endref

\ref \key M \by Mantovan, E. \paper \rm On certain unitary group Shimura varieties, Harvard Ph. D. thesis (2002)
\endref

\ref \key R \by Rapoport, M. \paper \rm Non-archimedean period
domains, Proceedings of the International Congress of
Mathematicians, Z\"urich, 1994, 423--434. (1995) \endref

\ref \key Rd \by Rodier, F. \paper \rm  Repr\'esentations de $GL(n,k)$ o\`u $k$ est un corps $p$-adique,
\linebreak {\it Ast\'erisque}, {\bf 92-93} (1982), 201-218; S\'em Bourbaki 1981-82, expos\'e no. 583 \endref

\ref \key RZ \by Rapoport, M. et T. Zink \paper \rm  {\it Period Spaces for $p$-divisible Groups}, Princeton:
Annals of Mathematics Studies {\bf 141} (1996) \endref

\ref \key Ro \by Rogawski, J. \paper \rm Representations of
$GL(n)$ and division algebras over a $p$-adic field, {\it Duke
Math. J.}, {\bf 50}, 161--196 (1983) \endref

\ref \key Sh \by Shahidi, F. \paper \rm  Local coefficients and
normalization of intertwining operators for $GL(n)$, {\it Comp.
Math.}, {\bf 48}, 271--295 (1983) \endref

\ref \key S \by Strauch, M. \paper \rm  On the Jacquet-Langlands correspondence in the cohomology of the
Lubin-Tate deformation tower, Preprintreihe SFB 478 (M\"unster), {\bf 72}, (1999) \endref

\ref \key SS \by Schneider, P. and U. Stuhler \paper \rm  The
cohomology of $p$-adic symmetric spaces, {\it Invent. Math.},
{\bf 105}, 47--122 (1991) \endref

\ref \key W \by Wedhorn, T. \paper \rm  The local Langlands correspondence for $GL(n)$ over $p$-adic fields,
lecture at the summer school on Automorphic Forms on GL(n) at the ICTP Trieste,
        ICTP Lecture Notes Series (to appear) \endref

\ref \key Z \by Zelevinsky,  A. V. \paper \rm  Induced
representations of reductive $p$-adic groups II:  on irreducible
representations of $GL(n)$, {\it Ann. Sci. E.N.S.}, {\bf 13},
165--210 (1980) \endref

\end